\documentclass[11pt]{amsart}
\usepackage{amsfonts,amsmath,amsthm,amssymb, color,amscd}
\usepackage{fullpage}

\numberwithin{equation}{section}

\newtheorem{thm}{Theorem}
\newtheorem*{thm*}{Theorem}

\newtheorem*{prop*}{Proposition}

\newtheorem{cor}[thm]{Corollary}

\newtheorem{lemma}[thm]{Lemma}

\newtheorem*{remark*}{Remark}
\newtheorem*{remarks*}{Remarks}

\newcommand{\ip}[1]{\langle #1 \rangle}

\newcommand{\restrictto}[2]{\left. #1 \right|_{#2}}
\newcommand{\ddtat}{\restrictto{\frac{d}{dt}}{t=0}}
\newcommand{\ddsat}{\restrictto{\frac{d}{ds}}{s=0}}
\newcommand{\partialddtat}{\restrictto{\frac{\partial}{\partial t}}{t=0}}

\begin{document}

\title[Homog. Ricci solitons]{Homogeneous Ricci solitons are algebraic}
\author[Michael Jablonski]{Michael Jablonski}
\thanks{ This work was supported in part by NSF grant DMS-1105647.}
\maketitle

\begin{abstract}
In this short note, we show that homogeneous Ricci solitons are algebraic.  As an application, we see that the generalized Alekseevskii conjecture is equivalent to the Alekseevskii conjecture.
\end{abstract}

\setcounter{section}{0}

\section{Introduction}

A Riemannian manifold $(M,g)$ is said to be a Ricci soliton if it satisfies the equation
	\begin{equation}\label{eqn: ricci soliton}   ric_g = cg + L_Xg\end{equation}
for some $c\in\mathbb R$ and some smooth vector field $X\in \mathfrak X(M)$.  Such metrics are of interest as they correspond to self-similar solutions of the Ricci flow 
	$$\frac{\partial}{\partial t} g = -2ric_g$$
That is, $g$ is the initial value of a solution to the Ricci flow of the form  $g_t = c(t) \varphi_t^*g$, where $c(t)\in \mathbb R$ and $\varphi_t \in \mathfrak{Diffeo}(M)$.  In this way, Ricci solitons are   geometric fixed points of the flow and so are  special  metrics.

Homogeneous Ricci solitons  arise naturally as limits under the Ricci flow 
 \cite{Lott:DimReductionAndLongTimeBehaviorOfRicciFlow,Lauret:RicciFlowForSimplyConnectedNilmanifolds}  and, independently, hold 
a distinguished place apart from other homogeneous metrics.  For example, nilmanifolds cannot admit Einstein metrics, but do often admit Ricci solitons \cite{Jensen:TheScalarCurvatureOfLeftInvariantRiemannianMetrics,Jablo:ModuliOfEinsteinAndNoneinstein}, Ricci solitons on nilmanifolds are precisely the minima of a natural geometric functional \cite{LauretNilsoliton},  and Ricci solitons are  metrics of maximal symmetry on certain solvmanifolds \cite{Jablo:ConceringExistenceOfEinstein}.

One natural kind of example arises as follows.  Consider a homogeneous space $G/K$ where $K$ is closed and connected.  For every derivation $D\in Der(\mathfrak g)$ such that $D:\mathfrak k \to \mathfrak k$, we have a well-defined map $D_{\mathfrak g/\mathfrak k} : \mathfrak g/\mathfrak k \to \mathfrak g / \mathfrak k$.  Denote such derivations of $\mathfrak g$ by $Der(\mathfrak g/\mathfrak k)$.  
A homogeneous Ricci soliton $(G/K,g)$ is called \emph{$G$-semi-algebraic} if the $(1,1)$ Ricci tensor is of the form
	\begin{equation}\label{eqn: definition of semi-algebraic soliton} Ric = cId + \frac{1}{2}( D_{\mathfrak g/\mathfrak k} + D_{\mathfrak g/\mathfrak k} {}^t) \end{equation}
on $\mathfrak g/\mathfrak k \simeq T_eG/K$, for some $c\in \mathbb R$ and some  $D\in Der(\mathfrak g/\mathfrak k)$.  
This definition is motivated by the idea of taking our family of diffeomorphisms $\{\varphi_t \}$ above to come from automorphisms of the group $G$ which leave $K$ invariant, see \cite{Jablo:HomogeneousRicciSolitons} or \cite{LauretLafuente:StructureOfHomogeneousRicciSolitonsAndTheAlekseevskiiConjecture}  for more details.

If our semi-algebraic Ricci soliton satisfies the seemingly stronger condition that $D_{\mathfrak g/\mathfrak k}$ is symmetric, then it is called a \emph{$G$-algebraic Ricci soliton}.  Up to this point, all known examples of semi-algebraic Ricci solitons were in fact algebraic and isometric to solvmanifolds.  (This follows from \cite{Jablo:HomogeneousRicciSolitons} together with \cite{LauretLafuente:OnhomogeneousRiccisolitons}.)  Further, it was known that every homogeneous Ricci soliton must be semi-algebraic relative to its full isometry group \cite{Jablo:HomogeneousRicciSolitons}.  We now present our main result.

\begin{thm}\label{thm: main theorem} Every $G$-semi-algebraic  Ricci soliton is necessarily $G$-algebraic.\end{thm}

\begin{cor} Let $(M,g)$ be a homogeneous Ricci soliton.  There exists a transitive group $G$, of isometries,  such that $M=G/K$ is a $G$-algebraic Ricci soliton.
\end{cor}

The theorem above resolves questions raised by Lafuente-Lauret \cite{LauretLafuente:StructureOfHomogeneousRicciSolitonsAndTheAlekseevskiiConjecture} and  He-Petersen-Wylie \cite{HePetersenWylie:WarpProdEinsteinMetricsOnHomogAndHomogRicciSolitons}.  
In these works, it was shown that one can always extend a simply-connected, algebraic soliton to an Einstein metric on a larger homogeneous space.  
There the goal was to relate the classical Alekseevskii conjecture on Einstein metrics to a more general version for Ricci solitons.  More precisely, they showed that (among simply-connected manifolds) the Alekseevkii conjecture for Einstein metrics is equivalent to the (apriori) more general conjecture in the case of 
algebraic Ricci solitons.  We state these conjectures for completeness. 

\begin{quote}
\textbf{Alekseevskii Conjecture:}  Every homogeneous  Einstein metric with negative scalar curvature is isometric to a simply-connected solvmanifold.
\end{quote}

\begin{quote} 
\textbf{Generalized Alekseevskii Conjecture:}  Every expanding  homogeneous  Ricci soliton is isometric to a simply-connected  solvmanifold.
\end{quote}

Until now, it was not clear if these conjectures were equivalent.  Applying \cite{LauretLafuente:StructureOfHomogeneousRicciSolitonsAndTheAlekseevskiiConjecture} or \cite{HePetersenWylie:WarpProdEinsteinMetricsOnHomogAndHomogRicciSolitons} in the simply-connected case together with \cite{Jablo:StronglySolvable} and the results here, we now know the following.

\begin{thm} The generalized Alekseevskii conjecture is equivalent to the Alekseevskii conjecture.\end{thm}

\begin{remark*}  It is important to note that the Alekseevskii conjecture stated above is a more modern, geometric version than that given in \cite{Besse:EinsteinMflds}.  The version given in \cite{Besse:EinsteinMflds} has the weaker, topological conclusion that a non-compact, homogeneous, Einstein space  is only diffeomorphic to $\mathbb R^n$.  It is still an open question as to whether  the classical version stated in \cite{Besse:EinsteinMflds} is equivalent to the stronger version we pose above.
\end{remark*}

\textit{Acknowledgments:}  It is our pleasure to thank Ramiro Lafuente for providing useful comments on a draft of this manuscript.

\section{Ricci solitons by type}

The analysis  of (homogeneous) Ricci solitons varies depending on which of the following   categories the metric falls into.  A Ricci soliton is called \emph{shrinking, steady, or expanding} (respectively) if the cosmological constant $c$ appearing in Eqn.~\ref{eqn: ricci soliton} satisfies $c>0$, $c=0$, or $c<0$ (respectively).

\subsection*{Shrinking solitons}  The simplest example of a non-Einstein, homogeneous, shrinker is obtained by considering a compact homogeneous Einstein space $M'$ (which necessarily has positive scalar curvature) and taking a product with $\mathbb R^n$, i.e. $M=M'\times \mathbb R^n$.  Here the vector field $X\in\mathfrak{X}(M)$ appearing in Eqn.~\ref{eqn: ricci soliton} generates a family of diffeomorphisms which simply dilate the $\mathbb R^n$ factor.  Examples of this type are called trivial Ricci solitons and a result of Petersen-Wylie \cite{Petersen-Wylie:OnGradientRicciSolitonsWithSymmetry} says that every homogeneous shrinking Ricci soliton is finitely covered by a trivial one.  Observe that such spaces are algebraic Ricci solitons.

\subsection*{Steady solitons}  A homogeneous steady soliton is necessarily flat.  This well-known fact is proved as follows.  Along the Ricci flow of any homogeneous manifold, the scalar curvature $sc$ evolves by  the ODE
	$$\frac{d}{d t}sc = 2 |Ric|^2$$
As the scalar curvature of a steady soliton  does not change along the flow, we see that the homogeneous, steady solitons are Ricci flat and so flat by   \cite{AlekseevskiiKimelfeld:StructureOfHomogRiemSpacesWithZeroRicciCurv}.  Such spaces are trivially algebraic Ricci solitons.

\subsection*{Expanding solitons}  Every homogeneous, expanding Ricci soliton is necessarily non-compact, non-gradient and all known examples of such spaces are isometric to solvable Lie groups with left-invariant metrics.  While there is no characterization in this case as nice as the previous two cases,  new structural results have  recently appeared in \cite{LauretLafuente:StructureOfHomogeneousRicciSolitonsAndTheAlekseevskiiConjecture}.  The results obtained there are essential in our proof and we briefly recall  those which we need.   

We first observe that it suffices to prove the theorem for simply-connected manifolds.  Now consider a simply-connected, expanding, semi-algebraic Ricci soliton on $G/K$.  As $G/K$ is endowed with a $G$-invariant metric, $Ad(K)$ is contained in a compact subgroup of $Aut(G)$ and so we have a decomposition $\mathfrak g = \mathfrak p \oplus \mathfrak k$, where $\mathfrak p$ is an $Ad(K)$-complement to $\mathfrak k$.  We fix the point $p = eK \in M = G/K$ and naturally  identify $\mathfrak p$ with $T_pM$ as follows
	$$X\in\mathfrak p  \quad \leftrightarrow \quad \ddsat exp(sX)\cdot p = \ddsat exp(sX)K.$$
Although there is more than one choice of $\mathfrak p$ that one can make, we apply the work \cite{LauretLafuente:StructureOfHomogeneousRicciSolitonsAndTheAlekseevskiiConjecture} in the sequel and so we choose, as they do, to have $B(\mathfrak k,\mathfrak p)=0$, where $B$ is the Killing form of $\mathfrak g$.

As $G/K$ admits an expanding Ricci soliton, we know from \cite{LauretLafuente:StructureOfHomogeneousRicciSolitonsAndTheAlekseevskiiConjecture} that the group $G$ decomposes as $N\rtimes U$ where $N$ is the nilradical and $U$ is a reductive subgroup which contains the stabilizer $K$.  
Thus the underlying manifold of $M$ may be considered as $N \times U/K$ and we naturally identify the point $p=eK\in G/K$ with $(e,eK)\in N \times U/K$.  The subalgebra $\mathfrak u$ contains a subspace $\mathfrak h$ which is 
complementary to $\mathfrak k$, and so we have $T_pM \simeq \mathfrak p = \mathfrak n \oplus \mathfrak h$.   Furthermore, $\mathfrak n$ and $\mathfrak h$ are orthogonal subspaces of $T_pM$.  For more details, see \cite{LauretLafuente:StructureOfHomogeneousRicciSolitonsAndTheAlekseevskiiConjecture}.

Denote the restriction of our metric $g$ to $\mathfrak p \simeq T_eG/K$  by $\ip{\cdot,\cdot}$.  Denote by $H\in\mathfrak p$ the `mean curvature vector' of $G/K$ defined by
	$$\ip{H,X} = tr\, (ad\, X) \quad \mbox{ for all } X\in\mathfrak p$$
Observe that  $H\in\mathfrak h$.  It is a useful fact that the subspace $\mathfrak h$ of $\mathfrak u$  is $(ad~H)$-stable \cite[Prop.~4.1]{LauretLafuente:StructureOfHomogeneousRicciSolitonsAndTheAlekseevskiiConjecture}.  If $D$ is the soliton derivation appearing Eqn.~\ref{eqn: definition of semi-algebraic soliton}, then we have 
	$$ D = -ad~H + D_1$$
where 
$D_1$ is the derivation which vanishes on $\mathfrak u$ and restricts to the nilsoliton derivation on $\mathfrak n$.  

In \cite[Prop.~4.14]{LauretLafuente:StructureOfHomogeneousRicciSolitonsAndTheAlekseevskiiConjecture},   several conditions are given for when a semi-algebraic Ricci soliton is actually algebraic.  One of those conditions is
	\begin{equation}\label{eqn: s(ad H)=0}    S(ad~H|_\mathfrak h) = 0\end{equation}
where $S(A) =\frac{1}{2}(A+A^t)$.  This is the technical result that we will prove, from which the theorem follows.

\section{The proof of theorem \ref{thm: main theorem}}

The soliton inner product $\ip{\cdot, \cdot }$ on  $T_pM$ above gives rise to a  natural  inner product on the endomorphisms of $T_pM$ given by $\ip{A,B} = tr(AB^t)$, where $B^t$ denotes the metric adjoint of $B$ relative to $\ip{\cdot,\cdot}$.

\begin{lemma} Using the above inner product on endomorphisms we have
	$$\ip{(0,ad~H|_\mathfrak h), Ric} = 0$$
where $(0,ad~H|_\mathfrak h)$ is the map on $T_pM$ defined as $0$ on $\mathfrak n$ and $ad~H|_\mathfrak h$ on $\mathfrak h$.
\end{lemma}

\begin{remark*} As has been observed by R.~Lafuente \cite{Lafuente:OnHomogeneousWarpedProductEinsteinMetrics}, our proof of the lemma holds more generally.  In fact, one simply needs the group to satisfy $G = U\ltimes N$ with $N$ nilpotent, $U$ reductive, and $K<U$, the metric to satisfy $N \perp U/K$ at $eK$, and the element $H$ may be replaced by any $Y\in\mathfrak u$ satisfying $[Y,\mathfrak k]\subset \mathfrak k$.
\end{remark*}

Before proving the lemma, we use it to  verify that Eqn.~\ref{eqn: s(ad H)=0} holds.

\begin{proof}[Verification of \ref{eqn: s(ad H)=0}]
Consider the mean curvature vector $H\in \mathfrak u$.  As $\mathfrak u$ is reductive, $ad~H|_\mathfrak u$ is traceless.  Furthermore, since $ad~H$ vanishes on the stabilizer $\mathfrak k$ (see Eqn.~26 of \cite{LauretLafuente:StructureOfHomogeneousRicciSolitonsAndTheAlekseevskiiConjecture}) and $\mathfrak u = \mathfrak k \oplus \mathfrak h$, we see that $tr~ad~H|_\mathfrak h = 0$.  Together with the above lemma we have
	\begin{eqnarray*}
	0 &=& \ip{(0,ad~H|_\mathfrak h), Ric}\\
		&=& \ip{(0,ad~H|_\mathfrak h), cId -S(ad~H) + D_1}\\
		&=& \ip{ad~H|_\mathfrak h , cId|_\mathfrak h - S(ad~H|_\mathfrak h)}\\
		&=& c\ tr(ad~H|_\mathfrak h) - tr~S (ad~H|_\mathfrak h)^2\\
		&=& 0 - tr~S(ad~H|_\mathfrak h)^2
	\end{eqnarray*}
Thus $S(ad~H|_\mathfrak h) = 0$, as claimed.
\end{proof}

We now prove the lemma by considering a certain deformation of the metric $g$ on $M$.  As $ad~H$ vanishes on $\mathfrak k$ and $K$ is connected, the family of automorphisms $\Phi_t = C_{exp(tH)} \in Aut(U)$ is the identity on $K$ and hence gives rise to  well-defined diffeomorphisms $\phi_t$ on $U/K$ given by
	$$\phi_t(uK) = \Phi_t(u)K \quad \mbox{ for } u\in U$$
Note that $(\Phi_t)_* =Ad(exp(tH)) =  e^{t\ ad~H} \in Aut(\mathfrak u)$.  On the manifold  $M = N\times U/K$, we consider the family of diffeomorphisms given by
	$$ \varphi_t  =(id, \phi_t) \quad \mbox{ on } N\times U/K$$
The deformations of $g$ of interest are $g_t = \varphi_t {}^* g$.

 As $\varphi_t$ fixes the point $p:= eK = (e,eK)\in M = N\times U/K$, and  scalar curvature is an invariant, we have
	$$\ddtat sc(\varphi_t {}^* g) _p = 0$$
We  use this in the following general equation which holds for any family of metrics $\{g_t\}$ with variation $h = \frac{\partial}{\partial t} g_t$ (see \cite[Lemma 3.7]{ChowKnopf})
	\begin{equation}\label{eqn: scalar curv variation}
	\frac{\partial}{\partial t} sc = - \Delta \overline H + div(div~h) -\ip{h,ric}
	\end{equation}
where in local coordinates we have 
	\begin{equation}\label{eqn: delta H}
	\Delta \overline{H} = g^{ij}g^{kl} \nabla_i \nabla_j h_{kl} 
	\end{equation}
and 
	\begin{equation}\label{eqn: div div h}
	div(div~h) = g^{ij}g^{kl} \nabla_i \nabla_k h_{jl}
	\end{equation}
Observe, at the point $p :=eK = (e,eK)$ of $M$ we have $\frac{\partial}{\partial t}|_{t=0}  (\varphi_t)_* = (0,ad~H|_\mathfrak h)$ and so the lemma follows from Eqn.~\ref{eqn: scalar curv variation} (evaluated at $p$) upon showing the terms $\Delta \overline{H}$ and $div(div\ h)$ vanish.

\begin{remark*}
Recall that, in local coordinates, we define the metric inverse $g^{ij}$ as the function satisfying  $\delta _i^l = g^{ij}g_{jl}$.  By choosing a frame which is $g$-orthonormal at every point, one would have that both  $g_{ij}$ and  $g^{ij}$ are the identity.  We make such a choice below.
\end{remark*}

To ease computational burden, we build a frame which is $g$-orthonormal at every point and exploits the property that our metric $g$ is $G$-invariant.  We start with an orthonormal basis of $T_pM$.  As $T_pM = \mathfrak n \oplus \mathfrak h$, we may choose a basis $\{e_i\}$ which is the union of an orthonormal basis of $\mathfrak n$ together with an orthonormal basis of $\mathfrak h$.

Next, we extend the basis $\{e_i\}$ to a local frame nearby to $p\in M$.  To do this, we first consider a slice $\mathfrak S$ of the right $K$ action on $G$ through $e\in G$.  That is, we have a submanifold $\mathfrak S$ of $G$ containing $e$ such that $\dim \mathfrak S = \dim G/K$ and the map
	$$ s\mapsto sK \quad  s\in \mathfrak S   $$
is a diffeomorphism of a neighborhood of $e\in \mathfrak S$ to a neighborhood of $eK \in G/K$.  Now, for $q\in M$ nearby to $p$, there exists $s\in \mathfrak S$ such that $q=s\cdot p$ and we define
	$$e_i(q) = s_* e_i,$$
where $s_*$ denotes the differential of the translation $s: p \mapsto q$.  We note that the frame is well-defined as our choice of $s\in\mathfrak S$ is unique, since $\mathfrak S$ is a slice.  Furthermore, our frame is $g$-orthonormal as $g$ is $G$-invariant.

Using the above choice of frame nearby to $p\in M$, we now  study Eqns.~\ref{eqn: delta H} and \ref{eqn: div div h}.  
We begin by computing the variation $h$ of $g_t = \varphi_t {}^* g$ in terms of 
$\{e_i\}$.  For a point $q\in M$ near $p$,
	\begin{equation}\label{eqn: general variation}
	h_{ij}(q) = \partialddtat (g_t)_{ij}(q) = \partialddtat (g_t) (e_i(q),e_j(q)) = \partialddtat  g( (\varphi_t)_* e_i(q), (\varphi_t)_*  e_j(q) ) 
	\end{equation}
Next we compute $(\varphi_t)_* v_q$ for a vector $v_q \in T_qM$.

As  $G=NU$, there exist $n\in N$ and $u\in U$ such that $s\in\mathfrak S$ may be written as $s=nu$ and $q = (nu)\cdot p$.  Furthermore, there exists $X\in \mathfrak p = \mathfrak n \oplus \mathfrak h$ such that  $v_q = (nu)_* \ddsat exp( s X) \cdot p$.   To understand Eqn.~\ref{eqn: general variation}, we analyze separately  the cases when $X$ is an element of $\mathfrak n$ or  of $\mathfrak h$.

For $X\in\mathfrak n$, we have
	\begin{eqnarray}
	(\varphi_t)_* v_q		&=& (\varphi_t)_* (nu)_* X \nonumber \\
									&=& \ddsat \varphi_t(  nu \ exp(sX)\cdot p  ) \nonumber \\
									&=& \ddsat \varphi_t( n \ u \ exp(sX) \ u^{-1}\ u  \cdot p  ) \nonumber  \\
									&=& \ddsat \varphi_t( n \ u \ exp(sX) \ u^{-1}, \ uK  )  \label{eqn: X in n}   \\
									&=& \ddsat ( n \ exp(sAd_uX), \Phi_t( u) K  ) \nonumber \\
									&=& \ddsat ( n \Phi_t(u)  \  \Phi_t(u)^{-1} \ exp(sAd_uX)  \Phi_t( u) K  ) \nonumber \\
									&=& \ddsat ( n \Phi_t(u) \  exp(sAd_{\Phi_t(u)^{-1}} Ad_u  X)  K  ) \nonumber \\
									&=& (n\Phi_t(u))_*  \ Ad_{\Phi_t(u)^{-1} u}  X  \nonumber
	\end{eqnarray}
Here we have used that $N$ is normal in $G$.  Note also  that $Ad_{\Phi_t(u)^{-1} u}  X \in \mathfrak n$. 

In the case when $X\in \mathfrak h \subset \mathfrak u$, we have
	\begin{eqnarray}
	(\varphi_t)_* v_q		&=& (\varphi_t)_* ( nu )_* X \nonumber \\
									&=& \ddsat \varphi_t(  nu \  exp(sX)\cdot p  ) \nonumber \\
									&=& \ddsat \varphi_t( n \ u \ exp(sX) \ K  ) \label{eqn: X in h}\\
									&=& \ddsat ( n \ \Phi_t(u \ exp(sX) )\ K  ) \nonumber \\
									&=& \ddsat ( n \ \Phi_t(u) \ exp(s (\Phi_t)_*X) )\ K  ) \nonumber \\
									&=& (n\Phi_t(u))_* (\Phi_t)_*X \nonumber
	\end{eqnarray}
Observe that since $ad~H$ preserves $\mathfrak h$ (\cite{LauretLafuente:StructureOfHomogeneousRicciSolitonsAndTheAlekseevskiiConjecture} Eqn. 32), $(\Phi_t)_*X \in \mathfrak h$ and so the the last line is consistent with our identification of $\mathfrak p = \mathfrak n \oplus \mathfrak h$ with $T_pM$.

From Eqns.~\ref{eqn: general variation}, \ref{eqn: X in n}, and \ref{eqn: X in h} we see that 
	\begin{enumerate}
	\item If $e_i \in \mathfrak n$ and $e_j \in \mathfrak h$, then $g_{ij}(q) = 0$.	
	\item If $e_i \in \mathfrak n$ and $e_j \in \mathfrak h$, then $h_{ij}(q) = 0$.
	\item If $e_i,e_j \in \mathfrak h$, then $h_{ij}(q)$ does not depend on $n$ and $u$, and so is constant in $q$.
	\item If $e_i,e_j \in \mathfrak n$, then $h_{ij}(q)$ does not depend on $n$, but does depend on $u$.
	\end{enumerate}
Using these observations, we see that the only possible non-zero terms of 
	$$div(div~h) = g^{ij}g^{kl} \nabla_i \nabla_k h_{jl}$$
are when $e_j,e_l \in \mathfrak n$ and $e_i, e_k \in \mathfrak h$.  However,  $(g_{\alpha \beta}) = Id$ implies $(g^{\alpha \beta})=Id$ and so $g^{kl}=0$.  This yields
	$$div(div~h) = 0$$
Next we study $\Delta \overline H = g^{ij}g^{kl} \nabla_i \nabla_j h_{kl} $.  
As above, the only possible non-zero terms occur when $e_k,e_l\in\mathfrak n$ and $e_i, e_j\in\mathfrak h$.  Further, as our frame is orthonormal, we have
	$$\Delta \overline H (q) = g^{ii} (q) g^{kk} (q)  ( \nabla_i \nabla_i h_{kk} )(q) = \sum_i  \left( \nabla_i \nabla_i \sum_k  h_{kk} \right) (q)$$
where the first sum is over the frame from $\mathfrak h$ and the second is over the frame from $\mathfrak n$.  
From Eqns.~\ref{eqn: general variation} and \ref{eqn: X in n} we have
	\begin{eqnarray*}
	h_{kk}(q) 			&=&   \partialddtat  g( (\varphi_t)_* e_k(q), (\varphi_t)_*  e_k(q) )  \\
							&=&   \partialddtat  \ip{ Ad_{\Phi_t(u)^{-1}  u}  ( e_k) , Ad_{\Phi_t(u)^{-1}  u}  ( e_k)   }  \\
							&=&   2 \  \ip{  e_k, (\ddtat Ad_{\Phi_t(u)^{-1}  u} ) ( e_k) } \\
							&=&   2 \  \ip{  e_k, ad~M ( e_k)}
	\end{eqnarray*}
where $M = \ddtat \Phi_t(u)^{-1}  u $.  To see that this last line makes sense, observe that $\Phi_t(u)^{-1}  u$ is a curve in $U$ with $\Phi_0(u)^{-1}u =e$ and thus $\ddtat \Phi_t(u)^{-1}  u \in \mathfrak u$.  

\begin{remark*}
Although $M$ is a function of $u$, we suppress this detail as it does not impact the rest of our proof.
\end{remark*}

We claim that $ad~M|_\mathfrak n$ is traceless.  To see this, we use that fact that $U$ being reductive and connected  implies $U = [U,U] Z(U)$, where $Z(U)$ is the center of $U$.  Thus, we may write $u=u_1u_2$ where $u_1\in [U,U]$ and $u_2\in Z(U)$.  As $u_2$ is central and $\Phi_t$ is an inner automorphism, $\Phi_t(u_2)=u_2$ and  
	$$\Phi_t(u)^{-1}u = \Phi_t(u_1)^{-1} u_1 \in [U,U]$$
This gives $ad~M \in ad~[\mathfrak u,\mathfrak u]$ from which our claim immediately follows.

Putting the above computations together, 
	\begin{eqnarray*}
	\Delta \overline H (q)	&=& \sum_i  \left( \nabla_i \nabla_i \sum_k  h_{kk} \right) (q) \\
										&=& 2 \sum_i \nabla_i \nabla_i \   tr\ ad\ M|_\mathfrak n \\
										&=& 0
	\end{eqnarray*}
which completes the proof of the lemma.

\bibliographystyle{amsalpha}

\end{document}